\input amstex
\documentstyle{amsppt}
\pageheight{43pc}
\pagewidth{28pc}
\NoBlackBoxes
\NoRunningHeads
\nologo
\TagsOnRight
\topmatter
\title{Two-Point Distortion Theorems for Harmonic Mappings}
\endtitle
\author{Martin Chuaqui, Peter Duren, and Brad Osgood}
\endauthor

\address Facultad de Matem\'aticas, P. Universidad Cat\'olica de Chile,
Casilla 306, Santiago 22, Chile
\endaddress
\email mchuaqui\@mat.puc.cl
\endemail

\address Department of Mathematics, University of Michigan, Ann Arbor,
Michigan 48109--1043
\endaddress
\email duren\@umich.edu
\endemail
 
\address Department of Electrical Engineering, Stanford University,
Stanford, California 94305
\endaddress
\email osgood\@ee.stanford.edu
\endemail

\thanks The authors are supported by Fondecyt Grant \# 1071019. 
\endthanks

\abstract 
In earlier work the authors have extended Nehari's well-known Schwarzian derivative 
criterion for univalence of analytic functions to a univalence criterion for canonical lifts 
of harmonic mappings to minimal surfaces.  The present paper develops some quantitative 
versions of that result in the form of two-point distortion theorems.   Along the way some 
distortion theorems for curves in ${\Bbb R}^n$ are given, thereby recasting a recent 
injectivity criterion of Chuaqui and Gevirtz in quantitative form.

\endabstract


\endtopmatter

\document

\flushpar
{\bf \S 1.  Introduction.}
\smallpagebreak

     The classical Koebe distortion theorem gives sharp bounds on the 
derivative of a normalized analytic univalent function.  Another measure of distortion 
is the distance $|f(z_1)-f(z_2)|$ between the images of two arbitrary points in 
the disk.   Some years ago, Blatter [3] gave a sharp lower bound for this distance 
in terms of the hyperbolic distance between $z_1$ and $z_2$.  More recently, Chuaqui 
and Pommerenke [10] found a sharp two-point distortion theorem for functions 
whose Schwarzian derivative satisfies Nehari's condition $|{\Cal S}f(z)|\leq2(1-|z|^2)^{-2}$.  
Their result may be viewed as a quantitative form of Nehari's univalence criterion.  The 
main purpose of the present paper is to carry out a similar analysis for harmonic mappings, 
or rather for their canonical lifts to minimal surfaces.  Along the way we obtain distortion 
theorems for curves in ${\Bbb R}^n$, thereby recasting an injectivity criterion of Chuaqui 
and Gevirtz [7] in quantitative form.  

     An important tool throughout the paper is the classical Sturm comparison theorem for 
solutions of linear differential equations of second order.  A good reference for this topic is 
the book of Birkhoff and Rota [2].

     The {\it Schwarzian derivative} of a locally univalent analytic
function is defined by
$$
{\Cal S}f=(f''/f')' - \tfrac12(f''/f')^2\,.
$$
It has the invariance property ${\Cal S}(T\circ f)={\Cal S}f$ for every
M\"obius transformation
$$
T(z)=\frac{az+b}{cz+d}\,, \qquad ad-bc\neq 0 \,.
$$
As a special case, ${\Cal S}(T)=0$ for every M\"obius transformation.  
A function $f$ has Schwarzian
${\Cal S}f=2\psi$ if and only if it has the form $f=u_1/u_2$ for some pair of 
independent solutions $w_1$ and $w_2$ of the linear differential
equation $w'' + \psi w=0$.  As a consequence, if  ${\Cal S}g={\Cal S}f$, then 
$g=T\circ f$ for some M\"obius transformation $T$.  In particular, 
M\"obius transformations are the only functions with  ${\Cal S}f=0$.

     In 1949, Nehari [14]  showed that if $f$ is analytic and locally univalent in the 
unit disk $\Bbb D$ and its Schwarzian satisfies either  $|{\Cal S}f(z)|\leq2(1-|z|^2)^{-2}$ 
or $|{\Cal S}f(z)|\leq{\pi}^2/2$ for all $z\in\Bbb D$, then $f$ is univalent 
in $\Bbb D$.  Pokornyi [16] then stated, and Nehari proved, that the 
condition $|{\Cal S}f(z)|\leq4(1-|z|^2)^{-1}$ also implies univalence.
Nehari [15] unified all three criteria by proving that $f$ is univalent 
under the general hypothesis $|{\Cal S}f(z)|\leq 2p(|z|)$, where $p(x)$ is a 
positive continuous even function defined on the interval $(-1,1)$, with the 
properties that $(1-x^2)^2p(x)$ is nonincreasing on the interval $[0,1)$ and no 
nontrivial solution $u$ of the differential equation $u''+pu=0$ has more than 
one zero in $(-1,1)$.  The last condition can be replaced by the equivalent 
requirement that some solution of the differential equation have no zeros in $(-1,1)$.    
We will refer to such functions $p(x)$ as {\it Nehari functions}.  

     It is clear from the Sturm comparison theorem that if $p(x)$ is a Nehari 
function, then so is  $cp(x)$ for any constant $c$ in the interval $0<c<1$.  
A Nehari function $p(x)$ is said to be {\it extremal} if $cp(x)$ is not a Nehari 
function for any constant $c>1$.  It was shown in [8] that some constant multiple of each 
Nehari function is an extremal Nehari function.  The functions $p(x)=(1-x^2)^{-2}$, 
$p(x)={\pi}^2/4$, and $p(x)=2(1-x^2)^{-1}$ are all extremal Nehari functions.  
Nonvanishing solutions of their corresponding differential equations are 
$u=\sqrt{1-x^2}$, $u=\cos(\pi x/2)$, and $u=1-x^2$, respectively.

     Ahlfors [1] introduced a notion of Schwarzian derivative 
for mappings of a real interval into ${\Bbb R}^n$, by formulating 
suitable analogues of the real and imaginary parts of ${\Cal S}f$ for 
analytic functions $f$.  A simple calculation shows that 
$$
\text{Re}\{{\Cal S}f\} = \frac{\text{Re}\{f'''\overline{f'}\}}{|f'|^2} 
-3 \,\frac{\text{Re}\{f''\overline{f'}\}^2}{|f'|^4} + \frac32 
\frac{|f''|^2}{|f'|^2}\,.
$$
For mappings $\varphi : (a,b)\mapsto {\Bbb R}^n$ of class $C^3$ with 
$\varphi'(x)\neq0$, Ahlfors defined the analogous expression
$$
S_1\varphi = \frac{\langle \varphi',\varphi'''\rangle}{|\varphi'|^2}
- 3\frac{\langle \varphi',\varphi''\rangle^2}{|\varphi'|^4} 
+ \frac32\frac{|\varphi''|^2}{|\varphi'|^2}\,, \tag1
$$
where $\langle\cdot\,,\cdot\rangle$ denotes the Euclidean inner product 
and $|{\bold x}|^2=\langle{\bold x},{\bold x}\rangle$ for ${\bold x}\in
{\Bbb R}^n$.   We will refer to $S_1\varphi$ as the {\it Ahlfors Schwarzian} of $\varphi$.  
As Ahlfors observed, it is invariant under postcomposition with M\"obius transformations; 
that is, under every composition of rotations, magnifications, translations, and inversions 
in ${\Bbb R}^n$.  

     In recent work, Chuaqui and Gevirtz [7] used the Ahlfors Schwarzian to give a criterion 
for injectivity of curves.  They proved the following theorem. 
\proclaim{Theorem A}  Let $p(x)$ be a continuous function such that the  
differential equation $u''(x)+p(x)u(x)=0$ admits no nontrivial solution 
$u(x)$ with more than one zero in $(-1,1)$. Let $\varphi : (-1,1)\mapsto 
{\Bbb R}^n$ be a curve of class $C^3$ with tangent vector 
$\varphi'(x)\neq 0$.  If $S_1\varphi(x)\leq2p(x)$, then $\varphi$ is injective.   
\endproclaim

     With the notation $v=|\varphi'|$, Chuaqui and Gevirtz also showed 
that 
$$
S_1\varphi = (v'/v)' - \tfrac12 (v'/v)^2 
+ \tfrac12 v^2 k^2 = {\Cal S}s + \tfrac12 v^2 k^2\,, \tag2
$$
where $s=s(x)$ is the arclength of the curve and $k$ is its scalar curvature, the 
magnitude of its curvature vector.

\bigpagebreak
\flushpar
{\bf \S 2.  Distortion of curves in ${\Bbb R}^n$.}
\smallpagebreak

     We now propose to give a sharpened form of Theorem A that expresses the 
injectivity in quantitative form by a two-point distortion inequality.   Closely related is an 
estimate for distortion in terms of the spherical derivative.  Here are our results.
\proclaim{Theorem 1} Let $p(x)$ be a positive continuous even function defined on the 
interval $(-1,1)$, with the property that no nontrivial solution $u$ of the differential equation 
$u''+pu=0$ has more than one zero in $(-1,1)$.   Let $F(x)$ be the solution to the differential 
equation ${\Cal S}F=2p$ determined by the conditions $F(0)=0$, $F'(0)=1$, and $F''(0)=0$.  
Let $\varphi : (-1,1)\mapsto {\Bbb R}^n$ be a curve of class $C^3$, normalized by 
$\varphi(0)=0$, $|\varphi'(0)|=1$, and $\langle\varphi'(0),\varphi''(0)\rangle=0$.  
If $S_1\varphi(x)\leq2p(x)$, then 
$$
\align
&\roman{(a)} \qquad |\varphi'(x)|\leq F'(x)\,, \qquad x\in(-1,1)\,,\qquad \text{and}\\
&\roman{(b)} \qquad \frac{|\varphi'(x)|}{1 + |\varphi(x)|^2} \leq \frac{F'(x)}{1 + F(x)^2}\,, \qquad 
x\in(-1,1)\,. 
\endalign
$$
\endproclaim
\proclaim{Theorem 2}  Let $p(x)$ and $F(x)$ be as in Theorem 1.  
If $\varphi : (-1,1)\mapsto {\Bbb R}^n$ is a curve of class $C^3$ with the property   
$S_1\varphi(x)\leq2p(x)$, then 
$$
\frac{|\varphi(x_1) - \varphi(x_2)|}{\{|\varphi'(x_1)||\varphi'(x_2)|\}^{1/2}}
\geq \frac{|F(x_1) - F(x_2)|}{\{F'(x_1)F'(x_2)\}^{1/2}}\,, \qquad x_1,x_2\in(-1,1)\,.
$$
\endproclaim

     The normalization required for the curve $\varphi$ in Theorem 1 can be achieved by 
postcomposing with a suitable M\"obius transformation.  Note that no such normalization is 
required for the two-point distortion result of Theorem 2.  Observe also that in Theorems 1 and 2 
it need not be assumed that $p(x)$ is a Nehari function.  In particular, no assumption is made  
that $(1-x^2)^2p(x)$ is nonincreasing, although this hypothesis will be essential in Theorem 3.

      Before passing to the proofs, it will be helpful to recall some properties of the function $F(x)$, 
which plays the role of extremal solution in Theorem A and in earlier work of Nehari.  
Since the function $p(x)$ of Theorem A is even, so is the solution $u_0$ 
of the differential equation $u''+pu=0$ with initial conditions $u_0(0)=1$ 
and $u_0'(0)=0$.   Therefore, $u_0(x)\neq0$ on $(-1,1)$, because otherwise 
it would have at least two zeros, contrary to hypothesis.  Thus the function 
$$
F(x) = \int_0^x 1/u_0(t)^2\,dt\,, \qquad -1 < x < 1\,, \tag3
$$
is well defined and satisfies the required initial conditions $F(0)=0$, $F'(0)=1$, and
$F''(0)=0$.  It also has the properties $F'(x)>0$ and $F(-x)=-F(x)$.  A calculation shows  
that $u_1=u_0F$ is an independent solution of $u''+pu=0$, and so $F=u_1/u_0$ has 
Schwarzian ${\Cal S}F=2p$.  Note also that $S_1F={\Cal S}F$, since $F$ is real-valued. 
In particular, $S_1F=2p$.  Finally, it should be noted that $F$ is strictly increasing on $(-1,1)$, 
because $F'(x)>0$.      

     For certain choices of $p(x)$ the function $F(x)$ can be calculated explicitly.  For 
instance, if $p(x)=(1-x^2)^{-2}$, then $u_0(x)=\sqrt{1-x^2}$ and so 
 $$
 F(x) = \int_0^x \frac{1}{1-t^2}\,dt = \frac12\log\frac{1+x}{1-x}\,.
 $$
 Similarly, for $p(x)={\pi}^2/4$ we have $u_0(x)=\cos(\pi x/2)$\,, so that
 $$
 F(x) = \int_0^x \sec^2(\pi t/2)\,dt = \frac{2}{\pi} \tan(\pi x/2)\,.
 $$
 If $p(x)=2(1-x^2)^{-1}$, then $u_0(x)=1-x^2$ and 
 $$
 F(x) = \int_0^x \frac{1}{(1-t^2)^2}\,dt = \frac14\,\log\frac{1+x}{1-x} + \frac12\,\frac{x}{1-x^2}\,.
 $$
 In such cases the distortion bounds in Theorems 1 and 2 take more concrete form.  
 For example, if $S_1\varphi(x)\leq\pi^2/2$, the inequality in Theorem 2 reduces to the 
 elegant form 
$$
\frac{|\varphi(x_1) - \varphi(x_2)|}{\{|\varphi'(x_1)||\varphi'(x_2)|\}^{1/2}} \geq
 \frac{2}{\pi} \,\sin\left(\frac{\pi}{2}|x_1 - x_2|\right)\,. 
$$
If $S_1\varphi(x)\leq2(1-x^2)^{-2}$, it says that 
$$
\frac{|\varphi(x_1) - \varphi(x_2)|}{\{|\varphi'(x_1)||\varphi'(x_2)|\}^{1/2}} \geq
\sqrt{(1-x_1^2)(1-x_2^2)}\ d(x_1,x_2)\,,
$$
where $d(x_1,x_2)$ is the hyperbolic distance between $x_1$ and $x_2$.
\demo{Proof of Theorem 1}  Part $(a)$ is in the paper by Chuaqui and Gevirtz [7] but 
we include the proof here for the sake of completeness.  It is known (and easy to verify) 
that if $g(x)$ is a real-valued function with $g'(x)>0$, then the function $u(x)=g'(x)^{-1/2}$ 
satisfies the differential equation $u''+\frac12({\Cal S}g)u=0$.  If we choose $g(x)=s(x)$, 
the arclength function along the given curve in ${\Bbb R}^n$, then $s'(x)=|\varphi'(x)|$ 
and $u(x)=|\varphi'(x)|^{-1/2}$ satisfies $u''+\frac12({\Cal S}s)u=0$.  Moreover, the 
normalization of the curve $\varphi$ implies that $u(0)=1$ and $u'(0)=0$.  But it follows 
from the relation (2) that ${\Cal S}s(x)\leq S_1\varphi(x)$, and by hypothesis 
$S_1\varphi(x)\leq2p(x)$, so we see that $\frac12{\Cal S}s(x)\leq p(x)$.  Thus it follows 
from the Sturm comparison theorem that $u(x)\geq u_0(x)$, which gives the 
inequality $(a)$.

     To prove $(b)$ we consider the inversion 
$$
\Phi(x) = \frac{\varphi(x)}{|\varphi(x)|^2}\,.
$$
Because the Ahlfors Schwarzian is M\"obius invariant, we see that $S_1\Phi=S_1\varphi$.  
On the other hand, we find as in the proof of Part (a) that the function $v(x)=|\Phi'(x)|^{-1/2}$ 
satisfies $v''+\frac12({\Cal S}s)v=0$, where now $s(x)$ denotes the arclength function along 
the curve $\Phi$, and 
$$
{\Cal S}s(x) \leq S_1\Phi(x) = S_1\varphi(x) \leq 2p(x)\,.
$$
A straightforward calculation shows that 
$$
|\Phi'(x)| = \frac{|\varphi'(x)|}{|\varphi(x)|^2}\,,
$$
so that $v(x)=|\varphi(x)|\,|\varphi'(x)|^{-1/2}$ and the normalization of the curve $\varphi$ 
implies that $v(0)=0$ and $v$ has a right-hand derivative  $v'(0)=1$.  On the other hand, the function $u_1=u_0F$ is a solution of $u''+pu=0$ with the same initial conditions $u_1(0)=0$ and $u_1'(0)=1$.  Therefore, the Sturm comparison theorem gives $v(x)\geq u_1(x)$ for $x>0$, or 
$$
\frac{|\varphi(x)|}{|\varphi'(x)|^{1/2}} \geq  \frac{|F(x)|}{F'(x)^{1/2}}  \tag4
$$
for $0\leq x<1$.  Since $v$ has a left-hand derivative $v'(0)=-1$, a similar argument shows 
that $-v(x)\geq u_1(x)$ for $x<0$, which implies that (4) holds also for $-1<x\leq0$.  Now 
square both sides of (4) and add the inequality of Part (a) in the form 
$1/|\varphi'(x)| \geq 1/F'(x)$ to obtain the desired result.  
\qed\enddemo

\demo{Proof of Theorem 2}   The proof is similar to that of Theorem 1.  Fixing any $x_1\in(-1,1)$, 
we now construct the inversion 
$$
\Phi(x) = \frac{\varphi(x) - \varphi(x_1)}{|\varphi(x) - \varphi(x_1)|^2}
$$ 
with respect to the point $\varphi(x_1)$.  By M\"obius invariance, $S_1\Phi=S_1\varphi$.  
The function $v(x)=|\Phi'(x)|^{-1/2}$ satisfies $v''+\frac12({\Cal S}s)v=0$, where $s(x)$ 
denotes the arclength function along the curve $\Phi$, and 
$$
{\Cal S}s(x) \leq S_1\Phi(x) = S_1\varphi(x) \leq 2p(x)\,.
$$
A calculation gives 
$$
|\Phi'(x)| = \frac{|\varphi'(x)|}{|\varphi(x) - \varphi(x_1)|^2}\,,
$$
so that  $v(x)=|\varphi(x) - \varphi(x_1)|\,|\varphi'(x)|^{-1/2}$.   Now $v(x_1)=0$ and a calculation 
shows that $v$ has right-hand derivative $v'(x_1)=|\varphi'(x_1)|^{1/2}$.  If $U(x)$ is the solution of the equation $u''+pu=0$ with $U(x_1)=0$ and $U'(x_1)=1$, the Sturm comparison theorem gives the 
inequality $|\varphi'(x_1)|^{-1/2}v(x)\geq U(x)$ for $x>x_1$.  To calculate the function $U(x)$, first let 
$$
H(x) = - \,\frac{1}{F(x)-F(x_1)}\,, \qquad \text{so that}\qquad H'(x) = \frac{F'(x)}{[F(x)-F(x_1)]^2}\,.
$$
Note that ${\Cal S}H= {\Cal S}F=2p$ by the M\"obius invariance of the Schwarzian.  
Thus by the general principle stated at the start of the proof of Theorem 1, the function 
$$
w(x) = H'(x)^{-1/2} = \frac{F(x) - F(x_1)}{F'(x)^{1/2}}
$$
satisfies the equation $w''+pw=0$ for $x>x_1$.  Also $w(x_1)=0$ and $w'(x_1)=F'(x_1)^{1/2}$.  This shows that $U(x)=F'(x_1)^{-1/2}w(x)$, so that the inequality $|\varphi'(x_1)|^{-1/2}v(x)\geq U(x)$ takes 
to the form 
$$
\frac{|\varphi(x) - \varphi(x_1)|}{\{|\varphi'(x_1)||\varphi'(x)|\}^{1/2}}
\geq \frac{|F(x) - F(x_1)|}{\{F'(x_1)F'(x)\}^{1/2}}\,, \qquad x_1 \leq x < 1\,.
$$
Now let $x=x_2$ to obtain the inequality of Theorem 2.
\qed\enddemo

       The bounds in Theorems 1 and 2 are sharp.  Equality occurs in all 
cases only when the curvature $\kappa=0$, so that the curve $\varphi$ is a straight line.  
Indeed, the relation (2) gives the inequality ${\Cal S}s(x) \leq S_1\Phi(x) $, with equality only 
when $\kappa=0$.   More precisely, in Theorem 1 equality occurs in either $(a)$ or $(b)$ at 
some point $x_0$ if and only if the portion of the curve $\varphi(x)$ between 0 and $\varphi(x_0)$ 
is a straight line that is parametrized so that $|\varphi'(x)|=F'(x)$ for all $x$ in the interval between 
0 and $x_0$.   In Theorem 2 equality occurs for a pair of points $x_1$ and $x_2$ if and only if the 
curve is a straight line between the points $\varphi(x_1)$ and $\varphi(x_1)$ that is parametrized 
so that  $|\varphi'(x)|=F'(x)$ for all $x$ in the interval between $x_1$ and $x_2$.

\bigpagebreak
\flushpar
{\bf \S 3.  Distortion of harmonic lifts.}
\smallpagebreak

     With the help of Theorem 2, we can now derive a two-point distortion inequality for the 
canonical lift of a harmonic mapping to a minimal surface.    A harmonic mapping is a 
complex-valued harmonic function $f(z)=u(z)+iv(z)$, for $z=x+iy$ in the unit disk $\Bbb D$ 
of the complex plane.   Such a mapping has a canonical decomposition $f=h+\overline{g}$, 
where $h$ and $g$ are analytic in $\Bbb D$ and $g(0)=0$.  

     According to the Weierstrass--Enneper formulas, a harmonic mapping $f=h+\overline{g}$ with 
$|h'(z)|+|g'(z)|\neq0$ lifts locally to a minimal surface described by conformal parameters if and only if 
its dilatation $\omega=g'/h'$ has the form $\omega=q^2$ for some meromorphic function $q$.  
The Cartesian coordinates $(U,V,W)$ of the surface are then given by
$$
U(z)=\text{Re}\{f(z)\}\,,\quad 
V(z)=\text{Im}\{f(z)\}\,,\quad
W(z)= 2\,\text{Im}\left\{\int_0^z 
h'(\zeta)q(\zeta)\,d\zeta\right\}\,.
$$
We use the notation $\widetilde{f}(z) = \bigl(U(z),V(z),W(z)\bigr)$ for the lifted mapping from 
$\Bbb D$ to the minimal surface.  The first fundamental form of the surface is 
$ds^2=\lambda^2|dz|^2$, where the conformal metric is $\lambda= |h'|+|g'|$\,.  The Gauss 
curvature of the surface at a point $\widetilde{f}(z)$ is 
$$
K = - \,\frac{1}{\lambda^2} \,\Delta (\log \lambda)\,,  
$$  
where $\Delta$ is the Laplacian operator.  Further information about 
harmonic mappings and their relation to minimal surfaces can be found 
in the book [12].
 
     For a harmonic mapping $f=h+\overline{g}$ with $\lambda(z)=|h'(z)|+|g'(z)|\neq0$, 
whose dilatation is the square of a meromorphic function, the {\it Schwarzian derivative} is 
defined [2] by the formula 
$$
{\Cal S}f = 2\bigl(\sigma_{zz} - {\sigma_z}^2\bigr)\,,  \qquad \sigma = \log\lambda\,,
$$
where  
$$
\sigma_z = \frac{\partial\sigma}{\partial z}  
= \frac12 \left(\frac{\partial\sigma}{\partial x} 
- i \frac{\partial\sigma}{\partial y}\right)\,, \qquad z = x+iy\,.
$$
If $f$ is analytic, it is easily verified that ${\Cal S}f$ reduces to the classical Schwarzian.  

     In our paper [5] we found the following criterion for the lift of a harmonic mapping to be 
univalent.

\proclaim{Theorem B}  Let $f=h+\overline{g}$ be a harmonic mapping of the 
unit disk, with $\lambda(z)=|h'(z)|+|g'(z)|\neq0$ and dilatation 
$g'/h'=q^2$ for some meromorphic function $q$.  Let $\widetilde{f}$ 
denote the Weierstrass--Enneper lift of $f$ to a minimal surface with 
Gauss curvature $K=K(\widetilde{f}(z))$ at the point $\widetilde{f}(z)$.  
Suppose that the inequality   
$$
|{\Cal S}f(z)| + {\lambda(z)}^2 |K(\widetilde{f}(z))| \leq 2p(|z|)\,, 
\qquad z\in\Bbb D\,,  \tag5
$$
holds for some Nehari function $p$. Then $\widetilde{f}$ is univalent in $\Bbb D$.
\endproclaim

     If $f$ is analytic, its associated minimal surface is the complex plane itself, with Gauss 
curvature $K=0$, and the result reduces to Nehari's theorem.  

     We can now sharpen Theorem B to express the univalence in quantitative form.  
Under the same hypotheses it turns out that the harmonic lift $\widetilde{f}$ 
actually satisfies a two-point distortion condition.  The inequality will involve the function 
$F$ determined by a Nehari function $p$ as in the formula (3).  In order to state 
the result in most elegant form, it will be convenient to assume that the given Nehari 
function is extremal, as defined in Section 1.  
\proclaim{Theorem 3}  Let $f$ be a harmonic mapping of the unit disk that has the properties 
specified in Theorem B, and let $\widetilde{f}$ be its canonical lift to a minimal surface.  
Suppose that the inequality \rom{(5)} holds for some extremal Nehari function $p$.  Then 
$\widetilde{f}$ satisfies the inequality
$$
|\widetilde{f}(z_1) - \widetilde{f}(z_2)| \geq \left\{\frac{\lambda(z_1)\lambda(z_2)}{F'(|z_1|)F'(|z_2|)}
\right\}^{1/2} d(z_1,z_2)\,, \qquad z_1,z_2\in\Bbb D\,,
$$
where $F(x)$ is defined by \rom{(3)} and $d(z_1,z_2)$ is the hyperbolic distance between the points 
$z_1$ and $z_2$.
\endproclaim
\demo{Proof}  The proof will apply Theorem 2.  The canonical lift  $\widetilde{f}$ onto a minimal surface 
$\Sigma$ defines a curve $\widetilde{f}: (-1,1) \to \Sigma\subset{\Bbb R}^3$.  As shown in [5], the Ahlfors Schwarzian of this curve satisfies the inequality 
$$
S_1\widetilde{f}(x) \leq |{\Cal S}f(x)| + {\lambda(z)}^2 |K(\widetilde{f}(x))|\,.
$$
Thus the hypothesis (5) tells us that $S_1\widetilde{f}(x) \leq 2p(x)$, and so by Theorem 2 we 
have the inequality 
$$
\frac{|\widetilde{f}(x_1) - \widetilde{f}(x_2)|}{\{\lambda(x_1)\lambda(x_2)\}^{1/2}}
\geq \frac{|F(x_1) - F(x_2)|}{\{F'(x_1)F'(x_2)\}^{1/2}}\,, \qquad x_1,x_2\in(-1,1)\,, \tag6
$$ 
since $|\widetilde{f}'(x)|=\lambda(x)$.  In order to extend the result to an arbitrary pair of distinct points $z_1, z_2\in\Bbb D$, we adapt a device due to Nehari [15].    It is here that the nonincreasing 
property of $(1-x^2)^2p(x)$ comes into play.  Suppose first that the hyperbolic geodesic $\gamma$ passing through $z_1$ and $z_2$ lies in the upper half-plane and is symmetric with respect to the imaginary axis.  Denote by $i\rho$ the midpoint of $\gamma$, so that $\rho>0$.   Then the M\"obius transformation 
$$
T(z) = \frac{i\rho -z}{1 + i\rho z}
$$
maps $\Bbb D$ onto itself and sends the segment $(-1,1)$ onto $\gamma$, with $T(x_1)=z_1$ 
and $T(x_2)=z_2$ for some pair of points $x_1$ and $x_2$.  The composite function 
$f_1(z)=f(T(z))$ is a harmonic mapping of the disk whose lift $\widetilde{f_1}=\widetilde{f}\circ T$ 
again maps $\Bbb D$ onto the minimal surface $\Sigma$.   Using the property of the Nehari 
function $p$ that $(1-x^2)^2p(x)$ is nonincreasing on $[0,1)$, we see as in [5] that (5) implies 
$$
|{\Cal S}f_1(x)| + {\lambda_1(x)}^2 |K(\widetilde{f_1}(x))| \leq 2p(x)\,, 
\qquad -1<x<1\,,  \tag7
$$
where $\lambda_1=|h_1'|+|g_1'|$ is the conformal factor associated with $f_1=h_1+\overline{g_1}$.  
It follows as before that $S_1\widetilde{f_1}(x))\leq 2p(x)$, and so by Theorem 2 the inequality 
(6) holds with $f$ replaced by $f_1$.  In other words,
$$
\frac{|\widetilde{f}(z_1) - \widetilde{f}(z_2)|}{\{\lambda(z_1)\lambda(z_2)\}^{1/2}}
\geq \frac{\{|T'(x_1)||T'(x_2)|\}^{1/2}|F(x_1) - F(x_2)|}{\{F'(x_1)F'(x_2)\}^{1/2}}\,. \tag8
$$ 

     We now develop a lower estimate for the right-hand side of the inequality (8) that depends 
explicitly on $z_1$ and $z_2$.   As shown in [9], the function $F$ coming from an 
extremal Nehari function $p$ has the property that $(1-x^2)F'(x)$ is nondecreasing on the 
interval $[0,1)$.  Since $F'$ is an even function with $F'(0)=1$, this shows that 
$(1-x^2)F'(x)\geq1$ on $(-1,1)$,  Therefore,
$$
|F(x_1) - F(x_2)| = \int_{x_1}^{x_2} F'(x)\,dx \geq  \int_{x_1}^{x_2}  \frac{1}{1-x^2}\,dx = d(x_1,x_2)\,.
$$
By M\"obius invariance of the hyperbolic metric, it follows that $|F(x_1) - F(x_2)|\geq d(z_1,z_2)$.  
On the other hand,  
$$
\frac{|T'(x)|}{1-|T(x)|^2} = \frac{1}{1-x^2} 
$$
and a simple calculation shows that $|T(x)|>|x|$, so that 
$$
(1-x_j^2) F'(x_j) \leq (1-|z_j|^2) F'(|z_j|) = (1-x_j^2) |T'(x_j)| F'(|z_j|)\,, \quad j=1,2\,.
$$
Consequently, 
$$
\frac{\{|T'(x_1)||T'(x_2)|\}^{1/2}|F(x_1) - F(x_2)|}{\{F'(x_1)F'(x_2)\}^{1/2}} \geq 
\frac{d(z_1,z_2)}{\{F'(|z_1|)F'(|z_2|)\}^{1/2}}\,, \tag9
$$
and the desired result follows in the special case where the geodesic $\gamma$ is 
symmetric with respect to the imaginary axis.  The general result now follows from the 
obvious fact that the right-hand side of (9) is invariant under rotation of the disk.  
This proves Theorem 3.  
\qed\enddemo
     It should be observed that the inequality is sharp for the Nehari function $p(x)=(1-x^2)^{-2}$, 
since $(1-x^2)F'(x)$ is constant in this case.  It may also be remarked that the restriction to 
extremal Nehari functions is not essential.   If $p$ is not extremal, then $p_1=cp$ is an extremal 
Nehari function for some constant $c>1$, and the inequality (5) holds {\it a fortiori} with $p$ 
replaced by $p_1$.  However, the function $F$ that occurs in the lower bound must be calculated 
in terms of $p_1$ rather than $p$.

\bigpagebreak
\flushpar
{\bf \S 4.  Distortion in the surface metric.}
\smallpagebreak

     Although Theorem 3 expresses the univalence of the harmonic lift $\widetilde{f}$ in 
quantitative form, its estimate of distortion does not lead to a covering theorem analogous 
to the classical Koebe one-quarter theorem (see for instance [11]).  For that purpose it is 
natural to replace the Euclidean metric by the surface metric 
$$
\rho(w_1,w_2) = \int_\Gamma ds  =  \int_\gamma \lambda(z)\,|dz|\,,
$$
where $\Gamma$ is a geodesic joining the points $w_1$ and $w_2$ on the minimal surface 
$\Sigma=\widetilde{f}(\Bbb D)$ and $\gamma={\widetilde{f}}^{-1}(\Gamma)$ is its preimage 
in the unit disk.  (More precisely, in case there is no such geodesic, $\rho(w_1,w_2)$ is defined 
as the infimum of the lengths of all curves joining the two points.)

     Here another extremal function comes into play, a companion of the function $F$ that 
enters into Theorem 3.   Given a Nehari function $p$, let $u_1$ be the solution of the differential equation $u''-pu=0$ with initial conditions $u_1(0)=1$ and $u_1'(0)=0$.   Since $p(x)>0$ and 
$u_1(0)>0$, the solution $u_1$ is convex and so $u_1(x)\geq1$ in $(-1,1)$.   Define 
$$
G(x) = \int_0^x 1/u_1(t)^2\,dt\,.
$$
Then, by the initial remark in the proof of Theorem 1, we see that  ${\Cal S}G=-2p$.  It is also 
clear that $G(0)=0$, $G'(0)=1$, and $G''(0)=0$.  With this notation, we are now prepared to state the distortion theorem.
\proclaim{Theorem 4} Let $f$ be a harmonic mapping of the unit disk that has the properties 
specified in Theorem B.  Let $\widetilde{f}$ be its canonical lift to a minimal surface 
$\Sigma=\widetilde{f}(\Bbb D)$, with conformal metric $\lambda$ and $ \sigma = \log\lambda$.  
Suppose in particular that $f$ satisfies the condition \rom{(5)} for some Nehari function $p$.  
Suppose further that $p(x)$ is nondecreasing on the interval $[0,1)$.  Then for $0<r<1$, 
$$
\min_{|z|=r} \rho\bigl( \widetilde{f}(z), \widetilde{f}(0)\bigr) \geq 
\frac{\lambda(0) G(r)}{1+|\sigma_z(0)|G(r)}\,,  \tag10
$$
In particular, the surface $\Sigma$ contains a metric disk of radius 
$$
R=  \frac{\lambda(0) G(1)}{1+|\sigma_z(0)|G(1)}
$$
centered at $\widetilde{f}(0)$.
\endproclaim

     Before embarking on the proof, we will examine the particular case where $p(x)=(1-x^2)^{-2}$.
Then ({\it cf.} [8]) it can be verified that 
$$
u_1(x) = \frac12 \sqrt{1-x^2}\left\{\left(\frac{1+x}{1-x}\right)^{\sqrt{2}/2} + 
\left(\frac{1-x}{1+x}\right)^{\sqrt{2}/2}\right\}
$$ 
and
$$
G(x) = \frac{1}{\sqrt{2}} \frac{(1+x)^{\sqrt{2}} - (1-x)^{\sqrt{2}}}{(1+x)^{\sqrt{2}} + (1-x)^{\sqrt{2}}}\,, 
\tag11
$$
with $G(1)=1/\sqrt{2}$.  In the classical case where $f(z)=z+a_2z^2+\dots$ is analytic and 
satisfies $|{\Cal S}f(z)|\leq 2(1-|z|^2)^{-2}$, the covering radius in Theorem 4 reduces to 
$$
R=\frac{1}{|a_2|+\sqrt{2}}\,.
$$
 But a result of Ess\'en and Keogh [13] gives the coefficient bound $|a_2|\leq\sqrt{2}$ in this 
 case, so we conclude from Theorem 4 that the image $f(\Bbb D)$ contains the disk 
 $|w|<\sqrt{2}/4$.  This estimate is sharp, as shown in [13], with extremal function 
 $$
 G^{\star}(z)=\frac{G(z)}{1+\sqrt{2}G(z)} = \frac{\sqrt{2}}{4} \,\left[1 - \left(\frac{1-z}{1+z}\right)^{\sqrt{2}}\,
 \right] = z - \sqrt{2}\,z^2 + \dots\,,
 $$ 
 which has Schwarzian ${\Cal S}G^{\star}(z)=-2(1-z^2)^{-2}$.   It was shown in [8] that $f(\Bbb D)$ contains the 
 larger disk $|w|<1/2$ if $|{\Cal S}f(z)|\leq 2(1-|z|^2)^{-2}$ and $a_2=0$.
 
      To prepare for a proof of Theorem 4, we now state a lemma that expresses the Ahlfors 
 Schwarzian of the lift to $\Sigma$ of a curve in the disk.  It is a slight generalization of a formula 
 in [5], where the underlying curve was taken to be the real interval $(-1,1)$.  The formula 
 also plays a role in [6], where the setting is different but the derivation is essentially the same.
 \proclaim{Lemma}  Let $\gamma(t)$ be an arclength-parametrized curve in $\Bbb D$ with 
 curvature $\kappa(t)$, and let $\varphi(t)=\widetilde{f}(\gamma(t))$ be its lift to a curve $\Gamma$ 
 on the surface $\Sigma=\widetilde{f}(\Bbb D)$.  Let $k_e(t)$ denote the normal component of 
 the curvature vector of \,$\Gamma$ with respect to $\Sigma$, and let $K(\varphi(t))$ be the 
 Gauss curvature of \,$\Sigma$ at the point $\varphi(t)$.  Then
 $$
 (S_1\varphi)(t) = \roman{Re}\left\{({\Cal S}f)(\gamma(t))\,{\gamma}'(t)^2\right\}  + \tfrac12 
 \lambda(\gamma(t))^2 \left[K(\varphi(t)) + k_e(t)^2\right] + \tfrac12 \kappa(t)^2\,. \tag12
 $$
 \endproclaim
 \demo{Proof of Theorem 4}  For fixed $r\in(0,1)$, let $z_0$ be a point on the circle $|z|=r$ 
 where the minimum distance $\rho\bigl( \widetilde{f}(z), \widetilde{f}(0)\bigr)$ is attained.  
 Then the geodesic $\Gamma$ that joins $\widetilde{f}(0)$ to $\widetilde{f}(z_0)$ lies on the 
 subsurface ${\Sigma}_r=\{\widetilde{f}(z) : |z|\leq r\}$.  Let $\gamma={\widetilde{f}}^{-1}(\Gamma)$ 
 be the preimage in ${\Bbb D}_r=\{z\in{\Bbb D} : |z|\leq r\}$, and let $L\geq r$ denote the arclength of 
 $\gamma$.  Let $\gamma(t)$ be the parametrization of $\gamma$ with respect to arclength, 
 with $\gamma(0)=0$,  and let $\varphi(t)= \widetilde{f}(\gamma(t))$ be the corresponding parametrization of $\Gamma$.  Finally let 
 $$
 v(t) = \lambda(\varphi(t)) = |{\varphi}'(t)|\,, \qquad \text{and let} \ \  s(t) = \int_0^t v(\tau)\,d\tau 
 $$
 denote the arclength along the curve $\Gamma$.  According to the relation (2), the Ahlfors Schwarzian 
 of $\varphi$ has the form 
 $$
 S_1\varphi = {\Cal S}s + \tfrac12 v^2 \bigl({k_i}^2 + {k_e}^2\bigr)\,, 
 $$
 where $k_i$ and $k_e$ denote respectively the tangential and normal components of curvature.  
 Comparing this with the expression (12) for $S_1\varphi$ given in the lemma, we conclude that 
 $$
  ({\Cal S}s)(t) =  \roman{Re}\left\{({\Cal S}f)(\gamma(t))\,{\gamma}'(t)^2\right\}  + \tfrac12 v(t))^2 
  |K(\varphi(t))| + \tfrac12 \kappa(t)^2\,,  \tag13
  $$
  since the tangential curvature $k_i$ vanishes along a geodesic.  But the univalence criterion (5) 
  implies that 
  $$
  \roman{Re}\left\{({\Cal S}f)(\gamma(t))\,{\gamma}'(t)^2\right\} \geq - |({\Cal S}f)(\gamma(t))| 
  \geq v(t)^2  |K(\varphi(t))| - 2p(|\gamma(t)|)\,. 
  $$
  Hence it follows from (13) that 
  $$
   ({\Cal S}s)(t) \geq \tfrac32 v(t)^2 |K(\varphi(t))| - 2 p(|\gamma(t)|) \geq - 2 p(|\gamma(t)|).
   $$ 
   Now observe that $|\gamma(t)|\leq t$ since $t$ is the arclength of the curve from $\gamma(0)=0$ 
   to $\gamma(t)$.  Therefore, $p(|\gamma(t)|)\leq p(t)$ because of the hypothesis that $p$ is 
   nondecreasing on the interval $[0,1)$, and we have proved that 
   $$
   ({\Cal S}s)(t) \geq - 2p(t)\,, \qquad 0\leq t\leq L_1=\min\{1,L\}\,. \tag14
   $$
   This is the inequality we will need for application of the Sturm comparison theorem. 
   
        For that purpose, first note that the function $w=v^{-1/2}$ is the solution of 
    $$
    w'' + \tfrac12  ({\Cal S}s) w = 0\,, \qquad w(0)=\lambda(0)^{-1/2}\,, \ \  w'(0)=-\tfrac12 v'(0)
   \,\lambda(0)^{-3/2}\,, 
   $$
   with
   $$
   w'(0) \leq |w'(0)| \leq |{\lambda}_z(0)| \,\lambda(0)^{-3/2}\,.
   $$
   Next consider the solution $u_2(t)$ of the differential equation 
   $$
   u'' - pu=0  \qquad \text{with} \ \ u_2(0)= \lambda(0)^{-1/2}\,, \ \  u_2'(0) =  
   |{\lambda}_z(0)| \,\lambda(0)^{-3/2}\,.
   $$
   Since $-p(t) \leq \tfrac12 ({\Cal S}s)(t) $ by (14), and also $u_2(0)=w(0)$ and $u_2'(0)\geq w'(0)$, it follows from the Sturm comparison theorem that 
   $$
   w(t)\leq u_2(t)\,, \qquad 0\leq t\leq L_1\,.
   $$

         Now let 
   $$
   H(x) = \int_0^x 1/u_2(t)^2\,dt\,, 
   $$
   and observe that ${\Cal S}H=-2p={\Cal S}G$, so that $H(x)=T(G(x))$ for some M\"obius 
   transformation $T$.  In order to calculate $T$ explicitly, note first that $T(0)=0$ since 
   $H(0)=G(0)=0$, so that $T$ has the form $T(x)=x/(ax+b)$ for some real parameters $a$ and $b$. 
   Writing 
   $$
   [aG(x) + b]\,H(x) = G(x)
   $$
   and differentiating, we find 
   $$
   aG'(x)H(x) +[aG(x) + b]\,H'(x) = G'(x)\,,
   $$
   so that 
   $$
   b = \frac{G'(0)}{H'(0)} = \frac{u_2(0)^2}{u_1(0)^2}= \frac{1}{\lambda(0)}\,.
   $$
   Another differentiation produces the relation $a = |\lambda_z(0)| \,\lambda(0)^{-2}$.
   This shows that 
   $$
   H(x) = \frac{\lambda(0)^2 G(x)}{|\lambda_z(0)|G(x) +  \lambda(0)} 
   = \frac{\lambda(0)G(x)}{1 + |{\sigma}_z(0)|G(x)}\,.   \tag15
   $$
   Consequently, for $0\leq t\leq L_1$ we have 
   $$
   \int_0^t v(\tau)\,d\tau =  \int_0^t w(\tau)^{-2}\,d\tau \geq  \int_0^t u_2(\tau)^{-2}\,d\tau = H(t)\,.
   $$
   Hence 
   $$
   \align
   \rho(\widetilde{f}(\zeta), \widetilde{f}(0)) &= \int_\gamma \lambda(z)\,|dz| = \int_0^L v(\tau)\,d\tau \\
   &\geq \int_0^r v(\tau)\,d\tau \geq H(r)\,.
   \endalign
   $$
   In view of the formula (15), this gives the inequality (12) stated in Theorem 4.
   \qed\enddemo
   
        The class of harmonic mappings considered in Theorem 4, satisfying in particular the 
   inequality (5), is invariant under precomposition $f\circ T$ with M\"obius self-mappings of 
   the disk when $p(x)=(1-x^2)^{-2}$.  This property yields an invariant formulation of Theorem 4, 
   virtually as a corollary.  
   \proclaim{Corollary} Let $f$, $\widetilde{f}$, $\lambda$, and $\sigma$ be as in Theorem 4, 
   and suppose that $p(x)=(1-x^2)^{-2}$, so that 
   $$
  |{\Cal S}f(z)| + {\lambda(z)}^2 |K(\widetilde{f}(z))| \leq \frac{2}{(1-|z|^2)^2}\,, \qquad z\in\Bbb D\,.
  $$
  Then for each fixed $\alpha\in{\Bbb D}$ and $0<r<1$, 
  $$
  \min_{\left|\frac{z-\alpha}{1-\overline{\alpha}z}\right|=r} 
  \rho\bigl( \widetilde{f}(z), \widetilde{f}(\alpha)\bigr) \geq 
  \frac{(1-|\alpha|^2)\lambda(\alpha) G(r)}{1+\left|(1-|\alpha|^2)\sigma_z(\alpha)
  -\overline{\alpha}\right|\,G(r)}\,,  \tag16
  $$
  where $G$ is defined by \rom{(11)}.
  \endproclaim
\demo{Proof}  Consider the harmonic mapping 
$$
f_1(z) = f(T(z))\,, \qquad \text{where} \ \ T(z)=\frac{z+\alpha}{1+\overline{\alpha}z}\,.
$$
Let $\widetilde{f_1}(z) = \widetilde{f}(T(z))$ be its harmonic lift, and let 
$$
\lambda_1(z) = \lambda(T(z))|T'(z)| = \lambda(T(z))\,\frac{1-|\alpha|^2}{|1 + \overline{\alpha}z|^2} 
$$
denote its conformal metric, with $\sigma_1=\log\lambda_1$.  Then $\lambda_1(0)=(1-|\alpha|^2)
\lambda(\alpha)$ and 
$$
{\sigma_1}_z(z) = \sigma_z(T(z))T'(z) - \frac{\overline{\alpha}}{1+ \overline{\alpha}z}\,,
$$
so that 
$$
{\sigma_1}_z(0) = (1-|\alpha|^2)\sigma_z(\alpha) - \overline{\alpha}\,.
$$
The inequality (16) now follows from (10) and the fact that the circles $|z|=r$ 
 and $\left|\frac{z-\alpha}{1-\overline{\alpha}z}\right|=r$ correspond under the mapping $T$.
 \qed\enddemo


\newpage

\Refs
\ref \no 1 \by L\. V\. Ahlfors \paper Cross-ratios and Schwarzian 
derivatives in ${\Bbb R}^n$ \inbook Complex Analysis: Articles dedicated 
to Albert Pfluger on the occasion of his 80th birthday (J\. Hersch and 
A\. Huber, editors)   
\publ Birkh\"auser Verlag, Basel  \yr 1988 \pages 1--15
\endref
\ref\no 2\by G\. Birkhoff and G\.-C\. Rota \book Ordinary Differential Equations \publ 
4th Edition, Wiley, New York \yr 1989
\endref
\ref\no 3\by C. Blatter \paper Ein Verzerrungssatz f\"ur schlichte Funktionen 
\jour Comment\. Math\. Helv\. \vol 53 \yr1978\pages651--659
\endref
\ref\no 4\by M\. Chuaqui, P\. Duren, and B\. Osgood \paper The Schwarzian
derivative for harmonic mappings \jour J\. Analyse Math\.\vol 91\yr2003\pages
329--351
\endref
\ref\no 5 \by M\. Chuaqui, P\. Duren, and B\. Osgood \paper Univalence 
criteria for lifts of harmonic mappings to minimal surfaces \jour J\. Geom\. Analysis 
\vol 17 \yr 2007 \pages 49--74
\endref
\ref\no 6 \by M\. Chuaqui, P\. Duren, and B\. Osgood \paper Injectivity criteria for holomorphic 
curves in ${\Bbb C}^n$ \jour Pure Appl\. Math\. Quarterly, to appear
\endref
\ref\no 7 \by  M\. Chuaqui and J\. Gevirtz \paper Simple curves in 
${\Bbb R}^n$ and Ahlfors' Schwarzian derivative \jour Proc\. Amer\. Math\. 
Soc\. \vol 132 \yr 2004 \pages 223--230
\endref
\ref\no 8 \by M\. Chuaqui and B\. Osgood \paper Sharp distortion theorems associated with the 
Schwarzian derivative  \jour J\. London Math\. Soc\. \vol 48 \yr 1993 \pages 289--298
\endref
\ref\no 9 \by M\. Chuaqui and B\. Osgood \paper Finding complete conformal 
metrics to extend conformal mappings \jour Indiana Univ\. Math\. J\. 
\vol 47 \yr 1998 \pages 1273--1292
\endref
\ref\no 10 \by M\. Chuaqui and Ch\. Pommerenke \paper Characteristic properties of Nehari 
functions \jour Pacific J\. Math\. \vol 188 \yr 1999 \pages 83--94
\endref
\ref\no 11 \by P\. L\. Duren \book Univalent Functions \publ
Springer--Verlag, New York  \yr 1983 
\endref
\ref\no 12 \by P\. Duren \book Harmonic Mappings in the Plane \publ
Cambridge University Press, Cambridge, U\. K\. \yr 2004 
\endref
\ref\no 13 \by M\. Ess\'en and F\. R\. Keogh \paper The Schwarzian derivative and estimates 
of functions analytic in the unit disc \jour Math\. Proc\. Cambridge Philos. Soc\. \vol 78 \yr 1975 
\pages 501--511
\endref
\ref\no 14 \by Z\. Nehari \paper The Schwarzian derivative and schlicht 
functions \jour Bull\. Amer\. Math\. Soc\. \vol 55 \yr 1949 \pages 
545--551
\endref
\ref\no 15 \by Z\. Nehari \paper Some criteria of univalence \jour Proc\. 
Amer\. Math\. Soc\. \vol 5 \yr 1954 \pages 700--704
\endref
\ref \no 16 \by V\. V\. Pokornyi \paper On some sufficient conditions 
for univalence \jour Dokl\. Akad\. Nauk SSSR \vol 79 \yr 1951 
\pages 743--746 (in Russian)
\endref

\endRefs
\enddocument